\documentclass{amsart}
\usepackage{geometry} 
\usepackage{url}
\usepackage{amsthm}
\usepackage{amssymb}
\usepackage{latexsym}
\usepackage{amsfonts}
\usepackage{amsmath}
\usepackage{amstext}
\usepackage{accents}
\usepackage{cancel}
\usepackage{eucal}
\usepackage{amscd}
\usepackage{hyperref}
\usepackage[shortalphabetic]{amsrefs}

\begin{document}

\title[Sixth order flow of curves]{A sixth order flow of plane curves with boundary conditions} 
\author[J. McCoy]{James McCoy*} \thanks{* Corresponding author}
\address{Institute for Mathematics and its Applications, University of Wollongong}
\email{jamesm@uow.edu.au}
\author[G. Wheeler]{Glen Wheeler}
\address{Institute for Mathematics and its Applications, University of Wollongong}
\email{glenw@uow.edu.au}
\author[Y. Wu]{Yuhan Wu}
\address{Institute for Mathematics and its Applications, University of Wollongong}
\email{yw120@uowmail.edu.au}

\thanks{The research of the first and second authors was supported by Discovery Project grant DP150100375 of the Australian Research Council.  The research of the third author was supported by a University of Wollongong Faculty of Engineering and Information Sciences Postgraduate research scholarship.}
\begin{abstract}
We show that small energy curves under a particular sixth order curvature flow
with generalised Neumann boundary conditions between parallel lines converge
exponentially in the $C^\infty$ topology in infinite time to straight lines.  
\end{abstract}

\maketitle

\section{Introduction} \label{S:intro}
\newtheorem{BCs}{Lemma}[section]
\newtheorem{main}[BCs]{Theorem}
\newtheorem{evlneqns}[BCs]{Lemma}

Higher order geometric evolution problems have received increasing attention in
the last few years. Particular geometric fourth order equations occur in
physical problems and enjoy some interesting applications in mathematics.
We mention in particular for curves the curve diffusion flow and $L^2$-gradient flow of the
elastic energy, and for surfaces the surface diffusion and Willmore flows.
Flows of higher even order than four have been less thoroughly investigated,
but motivation for them and their elliptic counterparts comes for example from
computer design, where higher order equations are desirable as they allow more
flexibility in terms of prescribing boundary conditions \cite{LX}.  Such
equations have also found applications in medical imaging \cite{UW}.

In this article we are interested in curves $\gamma$ meeting two parallel lines
with Neumann (together with other) boundary conditions evolving under the $L^2$
gradient flow for the energy
$$\int_\gamma k_s^2 ds \mbox{.}$$
Here $k_s$ denotes the first derivative of curvature with respect to the arc
length parameter $s$.  Particularly relevant to us is the corresponding
consideration of the curve diffusion and elastic flow in this setting in
\cite{WW}. Other relevant works on fourth order flow of curves with boundary
conditions are \cites{DLP14, DP14, L12}. Of course if one instead considers
closed curves without boundary evolving by higher order equations, these have
been more thoroughly studied; we mention in particular \cites{DKS02, EGBMWW14,
GI99, PW16, W13}.

The remainder of this article is organised as follows.  In Section
\ref{S:setup} we describe the set-up of our problem, the normal variation of
the energy and the boundary conditions.  We define our corresponding gradient
flow, discuss local existence and give the relevant evolution equations of
various geometric quantities.  We also state our main theorem in this part, Theorem \ref{T:main}.
In Section \ref{S:estimates} we state the relevant tools from analysis to be
used including an interpolation inequality valid in our setting.  Under the
small energy condition \eqref{E:SEC} below, we show that the winding number of
curves under our flow is constant and remains equal to zero.  We show further that
under this condition the length of the curve does not increase and the
curvature and all curvature derivatives in $L^2$ are bounded under the flow.
That these bounds are independent of time implies solutions exist for all time.
In Section \ref{S:exp} we show under a smaller energy assumption that in fact
the $L^2$ norm of the second derivative of curvature decays exponentially under
the flow.  As a corollary we obtain uniform pointwise exponential decay of
curvature and all curvature derivatives to zero.  A stability argument shows
that the solution converges to a unique horizontal line segment.  The
exponential convergence of the flow speed allows us to describe the bounded
region in which the solution remains under the flow.

\section{The Set-up} \label{S:setup}
Let $\gamma_0 : \left[ -1, 1\right] \rightarrow \mathbb{R}^2$ be a (suitably)
smooth embedded (or immersed) regular curve.  Denote by $ds$ the arc length
element and $k$ the (scalar) curvature.
We consider the energy functional
$$E\left[ \gamma\right] = \frac{1}{2} \int_\gamma k_s^2\, ds$$
where $k_s$ is the derivative of curvature with respect to arc length.  We are
interested in the $L^2$ gradient flow for curves of small initial energy with
Neumann boundary conditions.

Under the normal variation $\tilde \gamma = \gamma + \varepsilon F \nu$ a
straightforward calculation yields
\begin{equation} \label{E:1}
  \left. \frac{d}{d\varepsilon} E\left[ \tilde \gamma \right] \right|_{\varepsilon=0}
  = - 2 \int_\gamma \left( k_{s^4} + k^2 k_{ss} - \frac{1}{2} k \, k_s^2\right) F \, ds
  + 2 \left[ k_s F_{ss} + k_{ss} F_s + \left( k_{sss} + k^2 k_s \right) F \right]_{\partial \gamma} \mbox{.}
\end{equation}

`Natural boundary conditions' for the corresponding $L^2$-gradient flow would
ensure that the above boundary term is equal to zero.  However, this term is
rather complicated.  In view of the first term in \eqref{E:1}, we wish to take
\begin{equation} \label{E:speed}
  F=   \left( k_{s^4} + k^2 k_{ss} - \frac{1}{2} k \, k_s^2 \right)
  \end{equation}
and the corresponding gradient flow 
\begin{equation} \label{E:theflow}
  \frac{\partial \gamma}{\partial t} =   F \nu \mbox{.}
\end{equation}

Differentiating the Neumann boundary condition (see also \cite[Lemma 2.5]{WW} for example) implies 
\begin{equation} \label{E:2}
  0= -F_s\left( \pm 1, t\right) = -k_{s^5} - k k_s k_{ss} -k^2 k_{sss} + \frac{1}{2} k_s^3 \mbox{.}
\end{equation}

As in previous work, we will assume the `no curvature flux condition' at the boundary, 
\begin{equation} \label{E:noflux}
  k_s\left( \pm 1, t\right) =0 \mbox{.}
\end{equation}  

The boundary terms in \eqref{E:1} then disappear if we choose, for example, 
\begin{equation} \label{E:noflux2}
  k_{sss}\left( \pm 1, t\right)=0 \mbox{.}
  \end{equation}
  This is in a way a natural choice because equation \eqref{E:2} then implies $k_{s^5}\left( \pm1, t\right)=0$.  In fact by an induction argument we have

\begin{BCs}
With Neumann boundary conditions and also \eqref{E:noflux} and \eqref{E:noflux2} satisfied, a solution to the flow \eqref{E:theflow} satisfies $k_{s^{2\ell+1}}=0$ on the boundary for $\ell \in \mathbb{N}$.
\end{BCs}

Let us now state precisely the flow problem.  

Let $\eta_{\pm}\left( \mathbb{R}\right)$ denote two parallel vertical lines in
$\mathbb{R}^2$, with distance between them $\left| e\right|$.  We consider a
family of plane curves $\gamma: \left[ -1, 1\right] \times \left[ 0, T\right)
\rightarrow \mathbb{R}^2$ satisfying the evolution equation \eqref{E:theflow}
with normal speed given by \eqref{E:speed}, boundary conditions
$$\gamma\left( \pm 1, t\right) \in \eta_{\pm}\left( \mathbb{R} \right)$$
$$\left< \nu, \nu_{\eta_{\pm}} \right> \left( \pm 1, t\right) = 0$$
$$k_s\left( \pm 1, t\right) = k_{sss}\left( \pm 1, t\right) = 0$$
and initial condition
$$\gamma\left( \cdot, 0 \right) = \gamma_0\left( \cdot \right)$$
for initial smooth regular curve $\gamma_0$.

\begin{main} \label{T:main}
There exists a universal constant $C>0$ such that the following holds.
For the flow problem described above, if the initial curve $\gamma_0$ satisfies $\omega = 0$ and
\begin{equation} \label{E:SEC}
  \delta = 
           \left( \frac{\sqrt{1717} - 37}{174} \right) \pi^3
           - \left\| k_s \right\|_2^2 L_0^3
	   > 0\,,
\end{equation}
where $L_0$ is the length of $\gamma_0$, then the solution exists for all time
$T=\infty$ and converges exponentially to a horizontal line segment
$\gamma_\infty$ with dist$(\gamma_\infty,\gamma_0) < C/\delta$.
\end{main}

In the above statement and throughout the article we use $\omega$ to denote the \emph{winding number}, defined here as
\[
\omega := \frac{1}{2\pi} \int_{\gamma} k \, ds \mbox{.}
\]

\noindent \textbf{Remarks:} 
\begin{itemize}
  \item The condition \eqref{E:SEC} is not optimal.  We can relax it a little
	  but it might be possible to relax further.  Where our estimates hold
	  under a weaker energy assumption we will state them so. 
  \item The exponential decay facilitates an explicit estimate on the distance
	  of $\gamma_\infty$ to $\gamma_0$.
\end{itemize}

Local existence of a smooth regular curve solution $\gamma: \left[-1, 1\right]
\times \left[ 0, T\right) \rightarrow \mathbb{R}^2$ to the flow problem
$\gamma: \left[-1, 1\right] \times \left[ 0, T\right) \rightarrow \mathbb{R}^2$
is standard.  If $\gamma_0$ also satisfies compatibility conditions, then the
solution is smooth on $\left[ 0, T\right)$.  In this article we focus on the
case of smooth initial $\gamma_0$.  However, $\gamma_0$ may be much less
smooth; equation \eqref{E:theflow} is smoothing. We do not pursue this here.

Similarly as in \cite{WW} and elsewhere we may derive the following:

\begin{evlneqns} \label{T:evlneqns}
Under the flow \eqref{E:theflow} we have the following evolution equations
\begin{enumerate}
  \item[\textnormal{(i)}] $\frac{\partial}{\partial t} ds =  k F \,ds$;
  \item[\textnormal{(ii)}]  $\frac{\partial}{\partial t} k = F_{ss} + k^2 F$;
  \item[\textnormal{(iii)}]  $\frac{\partial}{\partial t} k_s =  F_{sss} + k^2 F_s + 3 k k_s F$;
  \item[\textnormal{(iv)}] 
  \begin{eqnarray*}
\partial_t k_{s^{l}}&=&k_{s^{l+6}}+\sum_{q+r+u=l}(c_{qru}^{1}k_{s^{q+4}}k_{s^{r}}k_{s^{u}}
+c_{qru}^{2}k_{s^{q+3}}k_{s^{r+1}}k_{s^{u}}\nonumber\\
&&+c_{qru}^{3}k_{s^{q+2}}k_{s^{r+2}}k_{s^{u}}
+c_{qru}^{4}k_{s^{q+2}}k_{s^{r+1}}k_{s^{u+1}})\nonumber\\
&&+\sum_{a+b+c+d+e=l}c_{abcde}k_{s^{a}}k_{s^{b}}k_{s^{c}}k_{s^{d}}k_{s^{e}}
\end{eqnarray*}
for constants $c_{qru}^{1}, c_{qru}^{2}, c_{qru}^{3}, c_{qru}^{4}, c_{abcde}\in \mathbb{R}$ with $a,b,c,d,e,q,r,u\geq 0.$
\end{enumerate}
\end{evlneqns}

\section{Controlling the geometry of the flow} \label{S:estimates}
\newtheorem{PSW}{Lemma}[section]
\newtheorem{interp}[PSW]{Proposition}
\newtheorem{WN}[PSW]{Lemma}
\newtheorem{Length}[PSW]{Lemma}
\newtheorem{curvature}[PSW]{Proposition}
\newtheorem{ksl}[PSW]{Proposition}

We begin with the following standard result for functions of one variable.
\begin{PSW}[Poincar\'{e}-Sobolev-Wirtinger (PSW) inequalities] \label{T:PSW}
 Suppose $f:\left[ 0, L \right] \rightarrow \mathbb{R}$, $L>0$ is absolutely continuous.  
 \begin{itemize}
   \item If $\int_0^L f \,ds =0$ then
 $$\int_0^L f^2 ds \leq \frac{L^2}{\pi^2} \int_0^L f_s^2 ds \mbox{ and } \left\| f\right\|_\infty^2 \leq \frac{2L}{\pi} \int_0^L f_s^2 ds \mbox{.}$$
  \item Alternatively, if $f\left( 0\right) = f\left( L \right) =0$ then
 $$\int_0^L f^2 ds \leq \frac{L^2}{\pi^2} \int_0^L f_s^2 ds \mbox{ and } \left\| f\right\|_\infty^2 \leq \frac{L}{\pi} \int_0^L f_s^2 ds \mbox{.}$$
 \end{itemize}
 \end{PSW}
 
 To state the interpolation inequality we will use, we first need to set up
 some notation.  For normal tensor fields $S$ and $T$ we denote by $S \star T$
 any linear combination of $S$ and $T$.  In our setting, $S$ and $T$ will be
 simply curvature $k$ or its arc length derivatives.  Denote by $P_n^m\left(
 k\right)$ any linear combination of terms of type $\partial_s^{i_1} k \star
 \partial_s^{i_2}k \star \ldots \star \partial_s^{i_n}k$ where $m=i_1 + \ldots+
 i_n$ is the total number of derivatives.
 
 The following interpolation inequality for closed curves appears in
 \cite{DKS02}, for our setting with boundary we refer to \cite{DP14}.
 
 \begin{interp} \label{T:int}
 Let $\gamma: I \rightarrow \mathbb{R}^2$ be a smooth closed curve.  Then for any term $P_n^m\left( k\right)$ with $n\geq 2$ that contains derivatives of $k$ of order at most $\ell-1$,
$$\int_I \left| P_\nu^\mu\left( k \right)\right| ds \leq c \, L^{1-m-n} \left\| k \right\|_2^{n-m} \left\| k \right\|_{\ell, 2}^{p}$$
where $p = \frac{1}{\ell} \left( m + \frac{1}{2} n - 1\right)$ and $c=c\left( \ell, m, n\right)$.  Moreover, if $m+ \frac{1}{2} < 2\ell+1$ then $p<2$ and for any $\varepsilon>0$,
\begin{equation*}
  \int_I  \left| P_n^m\left( k \right)\right| ds \leq \varepsilon \int_I \left| \partial_s^\ell k \right|^2 ds
   + c\, \varepsilon^{\frac{-p}{2-p}} \left( \int_I \left| k\right|^2 ds\right)^{\frac{n-p}{2-p}} + c\left( \int_I \left| k\right|^2 ds \right)^{m+n-1}\mbox{.}
\end{equation*}
\end{interp}

Our first result concerns the winding number of the evolving curve $\gamma$.
In view of the Neumann boundary condition, in our setting the winding number
must be a multiple of $\frac{1}{2}$.

\begin{WN} \label{T:WN}
	Under the flow \eqref{E:theflow}, $\omega(t) = \omega(0)$.
\end{WN}

\noindent \textbf{Proof:} We compute using Lemma \ref{T:evlneqns} (i)
$$\frac{d}{dt} \int k\, ds = -\int F_{ss} ds - \int k^2 F ds + \int k^2 F ds =0 \mbox{,}$$
so $\omega$ is constant under the flow.\hspace*{\fill}$\Box$
~\\


\newcommand{\vn}[1]{\lVert #1 \rVert}

\noindent \textbf{Remarks:}
\begin{itemize}
\item It follows immediately that the average curvature $\overline{k}$ satisfies
$$\overline{k} := \frac{1}{L} \int_\gamma k\, ds \equiv 0$$
under the flow \eqref{E:theflow}.  This is important for applying the inequalities of Lemma \ref{T:PSW}.
\item Unlike the situation in \cite{WW}, here small energy does not automatically imply that the winding number is close to zero.
	Indeed, one may add loops (or half-loops) of circles that contribute an arbitrarily small amount of the energy $L^3\vn{k_s}_2^2$.
	Note that such loops must all be similarly oriented, as a change in
	contribution from positive to negative winding will necessitate a
	quantum of energy (for example a figure-8 style configuration with
	$\omega=0$ can not have small energy despite comprising essentially of
	only mollified arcs of circles).
\end{itemize}

Next we give an estimate on the length of the evolving curve in the case of
small initial energy.  Of course, this result does not require the energy as
small as \eqref{E:SEC}.

\begin{Length} \label{T:L}
	Under the flow \eqref{E:theflow} with $\omega(0) = 0$,
$$\frac{d}{dt} L \left[ \gamma\left( t\right) \right] \leq 0 \mbox{.}$$
\end{Length}

\noindent \textbf{Proof:} We compute using integration by parts
\begin{equation*}
\frac{d}{dt} L \left[ \gamma\left( t\right) \right] = - \int k F\, ds = - \int k_{ss}^2 ds + \frac{7}{2} \int k^2 k_s^2 ds
\leq -\left[ 1 - \frac{7 L^3}{\pi^3} \left\| k_s\right\|_2^2 \right] \int k_{ss}^2 ds 
\end{equation*}
 where we have used Lemma \ref{T:PSW}.  The result follows by the small energy assumption.\hspace*{\fill}$\Box$

~\\
Thus under the small energy assumption we have the length of the evolving curve bounded above and below:
$$\left| e\right| \leq L\left[ \gamma\right] \leq L_0 \mbox{.}$$

We are now ready to show that the $L^2$-norm of curvature remains bounded, independent of time.

\begin{curvature} \label{T:k}
	Under the flow \eqref{E:theflow} with $\omega(0) = 0$, there exists a universal $C>0$ such that
$$\left\| k\right\|_2^2 \leq \left. \left| k \right\|_2^2 \right|_{t=0} + C \mbox{.}$$
\end{curvature}

\noindent \textbf{Proof:} Using integration by parts, Lemma \ref{T:evlneqns} and the interpolation inequality Proposition \ref{T:int}
\begin{multline*}
  \frac{d}{dt} \int k^2 ds = - 2\int k_{s^3}^2 ds + 5 \int k_{ss}^2 k^2 ds + 5 \int k_{ss} k_s^2 k \, ds
 + \int k_{ss} k^5 ds - \frac{1}{2} \int k_s^2 k^4 ds\\
  \leq \left( -2 + 3 \varepsilon\right) \int k_{s^3}^2 ds + C \left\| k \right\|_2^{14}
  \leq -\frac{\pi^6}{L_0^6} \int k^2 ds + \frac{C\pi^7}{\left| e\right|} \mbox{,}
  \end{multline*}
  where we have also used Lemma \ref{T:PSW} and the length bounds.  The result follows.\hspace*{\fill}$\Box$
 ~\\

Moreover, we may show similarly using the evolution equation for $k_{s^\ell}$ that all derivatives of curvature are bounded in $L^2$ independent of time.

\begin{ksl} \label{T:ksl}
	Under the flow \eqref{E:theflow} with $\omega(0) = 0$, there exists a universal $C>0$ such that, for all $\ell \in \mathbb{N}$,
$$\left\| k_{s^\ell} \right\|_2^2 \leq \left. \left| k_{s^\ell} \right\|_2^2 \right|_{t=0} + C \mbox{.}$$
\end{ksl}

Pointwise bounds on all derivatives of curvature follow from Lemma \ref{T:PSW}.
It follows that the solution of the flow remains smooth up to and including the final time, from which we may (if $T<\infty$) apply again local existence.
This shows that the flow exists for all time, that is, $T=\infty$.

\section{Exponential convergence} \label{S:exp}
\newtheorem{kssevln}{Lemma}[section]
\newtheorem{kssexp}[kssevln]{Corollary}

\begin{kssevln}\label{T:kssevln}
	Under the flow \eqref{E:theflow} with $\omega(0) = 0$,
$$\frac{d}{dt} \int k_{ss}^2 ds \leq \left[ -2 + \frac{74L^3}{\pi^3} \left\| k_s\right\|_2^2 + \frac{174 L^6}{\pi^6} \left\| k_s \right\|_2^4\right] \left\| k_{s^5}\right\|_2^2 - \frac{3}{L} \left\| k_{ss}\right\|_2^4 \mbox{.}$$
\end{kssevln}

Under the small energy condition of Theorem \ref{T:main}, the coefficient of $\left\| k_{s^5}\right\|_2^2$ above is bounded above by $-\delta$.  Using also Lemma \ref{T:PSW} we obtain the following.

\begin{kssexp} \label{T:kssexp}
Under the flow, there exists $\delta >0$ such that
$$\frac{d}{dt} \left\| k_{ss}\right\|_2^2 \leq - \delta \left\| k_{ss}\right\|_2^2  \mbox{.}$$
It follows that $\left\| k_{ss}\right\|_2^2$ decays exponentially to zero.
\end{kssexp}

\noindent \textbf{Remark:}
The small energy condition of Theorem \ref{T:main} may be slightly relaxed by
estimating the last term on the right hand side of Lemma \ref{T:kssevln} using
the H\"{o}lder inequality.\\

\noindent \textbf{Completion of the proof of Theorem \ref{T:main}:} Exponential
decay of $\left\| k_{ss}\right\|_2^2$ implies exponential decay of $\left\| k
\right\|_2^2$, $\left\| k_s \right\|_2^2$, $\left\| k\right\|_\infty$, $\left\|
k_s\right\|_\infty$ via Lemma \ref{T:PSW}.  Exponential decay of $\left\|
k_{s^\ell} \right\|_2$ and $\left\| k_{s^\ell} \right\|_\infty$ then follows by
a standard induction argument involving integration by parts and the curvature
bounds of Propositions \ref{T:k} and \ref{T:ksl}.  That $\left\| k_s
\right\|_2^2 \rightarrow 0$ implies subsequential convergence to straight line
segments (horizontal, in view of boundary conditions).  A stability argument
(see \cite{WW} for the details of a similar argument) gives that in fact the
limiting straight line is unique; all eigenvalues of the linearised operator
    $$\mathcal{L} u = u_{x^6}$$
are negative apart from the first zero eigenvalue, which corresponds precisely
to vertical translations. By Hale-Raugel's convergence theorem \cite{HR92} 
uniqueness of the limit follows.
Although we don't know the precise height of the
limiting straight line segment, we can estimate a-priori its distance
from the initial curve, since
\begin{equation*}
  \left| \gamma\left( x, t\right) - \gamma\left( x, 0\right) \right| = \left| \int_0^t \frac{\partial \gamma}{\partial t} \left( x, \tau\right) d \tau \right| \\
   \leq \int_0^t \left| F \right| d\tau \leq \frac{C}{\delta} \left( 1 - e^{-\delta t}\right) \mbox{.}
\end{equation*}

\end{document}